\documentclass{article}
\title{On nearly K\"ahler geometry 
}
\author{Paul-Andi Nagy}
\date{\today}
\usepackage{amsfonts}
\usepackage{amssymb}
\newtheorem{teo}{Theorem}[section]
\newtheorem{lema}{Lemma}[section]
\newtheorem{pro}{Proposition}[section]

\newtheorem{rema}{Remark}[section]
\newtheorem{coro}{Corollary}[section]
\newtheorem{nr}{}[section]
\begin{document}
\maketitle
\abstract{\normalsize We consider complete nearly K\"ahler manifolds with the canonical 
Hermitian connection. We prove some metric properties of strict nearly K\"ahler 
manifolds and give a sufficient condition for the reducibility of the canonical 
Hermitian connection. A holonomic condition for a nearly K\"ahler manifold to be 
a twistor space over a quaternionic K\"ahler manifold is given. This enables us to give 
classification results in 10-dimensions.}
\large
\section{Introduction}
Nearly K\"ahler (briefly NK) geometry is related to the concept of weak holonomy, 
introduced by A. Gray \cite{Gray2} in 1971. He proved that among those groups 
acting transitively on the sphere there are only 3 groups, namely
$$ U(n) \ \mbox{in dimension} \ 2n, \ G_2 \ \mbox{in dimension} \ 7, \ Spin(9) \ \mbox{in dimension} \ 16 $$
that can occur as weak holonomy groups and produce 
other geometries than the classical holonomy approach. Nearly K\"ahler 
geometry corresponds to weak holonomy $U(n)$ and was intensively studied in the seventies 
by Gray \cite{Gray4,Gray1}. Also note that the class of NK-manifolds appears naturally as one 
of the sixteen classes of almost Hermitian manifolds described by the Gray-Hervella 
classification \cite{Gray3}.\par
Recent interest for the study of such manifolds can be justified by the fact that in dimension $6$ 
nearly K\"ahler manifolds are related to the existence of a Killing spinor (see 
\cite{Grun1}). 
Furthermore, nearly K\"ahler manifolds provide a natural example 
of almost Hermitian manifolds admitting a Hermitian connection with totally skew symmetric 
torsion. From this point of view they are of interest in string theory 
(see \cite{Friedrich2}). \par
The aim of this paper is to investigate a number of properties of NK-manifolds related 
to the reducibility of the canonical Hermitian connection. We begin by 
proving a decomposition result which allows us to restrict our attention to 
strict NK-manifolds (see section 1). Our first main result 
is the following.
\begin{teo} Let $(M^{2n},g,J)$ a complete, strict nearly K\"ahler manifold. Then 
the following hold : \\
(i) If $g$ is 
not an Einstein metric then the canonical Hermitian connection has reduced holonomy. \\
(ii) The metric $g$ has positive Ricci curvature, hence $M$ is compact with finite
fundamental group. \\
(iii) The scalar curvature of the metric $g$ is a strictly positive constant.
\end{teo}
The previous theorem is a synthesis of the results contained in section 2. 
\par 
Let us recall now that one main class of examples of NK manifolds is formed by the so called 
$3$-symmetric spaces \cite{Gray4}. 
Other examples are provided by total spaces of Riemannian submersions with totally geodesic 
fibers admitting a compatible K\"ahler structure. These manifolds admit a canonical NK structure such that 
the canonical Hermitian connection has reduced holonomy (see section 3). In particular 
twistor spaces over positive quaternion-K\"ahler manifolds (here positive means of positive scalar curvature) 
have canonical NK-structures, a result already proven in \cite{Al1}. See also \cite{Mu} for the case of twistor bundles 
over $4$-manifolds. \par
In the second part of this paper we are concerned with the the study of the most simple case of 
reducible NK-geometry which is the following :
\begin{teo}
Let $(M^{2n}, g, J)$ be a complete, strict nearly K\"ahler manifold. If the holonomy group 
of the canonical Hermitian connection is contained in 
$U(1) \times U(n-1)$ then $M$ is the twistor space of a positive quaternionic-K\"ahler manifold endowed 
with its canonical NK-structure.
\end{teo}
In $6$-dimensions, the theorem 1.2 was already proven by a different 
method in  \cite{Moro1}. Our approach consists in showing that the torsion of the canonical 
Hermitian connection has to be of special algebraic type with respect to the holonomy 
decomposition. This will be done 
in section 4. Then, using standard arguments one can show that $M$ carries 
a complex contact structure and a K\"ahler-Einstein metric. The conclusion follows 
by a theorem of LeBrun (see section 5). \par
As a corollary of theorem 1.2 
we obtain a structure result in $10$-dimensions. Note that in $8$-dimensions it 
was already known by Gray \cite{Gray1} that there are no strict NK-manifolds.
\begin{coro}
Let $(M^{10},g,J)$ be a complete NK-manifold. Then either the universal cover of $M$ 
is a Riemannian product of a K\"ahler surface with a six dimensional NK-manifold, 
either $M$ is the twistor space of a positive, $8$-dimensional quaternionic K\"ahler manifold equipped 
with its canonical NK structure.
\end{coro}
Using results from \cite{Sal1} (see also \cite{Sal2}) we know that the only 
positive quaternionic-K\"ahler manifolds of $8$-dimensions are the symmetric spaces $P\mathbb{H}^2, 
\mathbb{G}r_{2}(\mathbb{C}^4), G_2 \slash SO(4)$ with their 
canonical metrics. Hence their twistor spaces, which are described in \cite{Sal2}, equipped with 
the canonical NK structure exhaust the list of complete, strict $NK$-manifolds of dimension $10$.
\section{Nearly K\"ahler geometry}
A nearly K\"ahler manifold is an almost hermitian manifold $(M^{2n},g,J)$ such 
that 
$$ (\nabla_X J)X=0$$
for every vector field $X$ on $M$ (here $\nabla$ denotes the Levi-Civita 
connection associated to the metric $g$). A NK manifold is called {\it{strict}} if $\nabla_XJ \neq 0$ for every 
$X \in TM, X \neq 0$. \par
Recall that the tensor $\nabla J$ has a number of important 
algebraic properties that can be summarized as follows : the tensors $A$ and $B$ defined for 
$X,Y,Z$ in $TM$ by $A(X,Y,Z)=<(\nabla_XJ)Y,Z>$ and $B(X,Y,Z)=<(\nabla_XJ)Y,JZ>$ are skew-symmetric and 
have type $(0,3)+(3,0)$ as real $3$-forms.
Denote by $Ric$ the Ricci tensor of the metric $g$ and by $Ric^{\star}$ its star 
version, that is the operator defined by 
$$<Ric^{\star}(X),Y>=\frac{1}{2}\sum \limits_{i=1}^{2n}R(X,JY,e_i,Je_i)$$
where $R$ is the curvature tensor of $(M,g)$ and $\{e_1,\ldots, e_{2n}\}$ a local frame 
field. The difference of these tensors, to be denoted by $r$, is given by the formula (see 
\cite{Gray1}) : 
$$ <rX,Y>=\sum \limits_{i=1}^{2n} <(\nabla_{e_i}J)X, (\nabla_{e_i}J)Y>.$$
Obviously $r$ is symmetric, positive and commutes with $J$. 
Another object of particular importance is the canonical hermitian connection defined by 
$$ \overline{\nabla}_XY=\nabla_XY+\frac{1}{2} (\nabla_XJ)JY.$$
It is easy to see that $\overline{\nabla}$ is the unique Hermitian con nection on $M$ with totally skew-symmetric 
torsion (see for example \cite{Friedrich2}). Note that the torsion of $\overline{\nabla}$ 
given by $T(X,Y)=(\nabla_XJ)JY$ vanishes iff $(M,g,J)$ is a K\"ahler manifold. \par 
The tensor $r$ has strong geometric properties. To begin, we have : 
\begin{nr} \hfill
$\overline{\nabla}r=0. \hfill $
\end{nr}
In fact, A. Gray proved in \cite{Gray1} that for all $X,Y,Z$ in $TM$ we have 
$$2<(\nabla_Xr)Y,Z>=\\<r(\nabla_XJ)Y,JZ)>+<r(JY),(\nabla_XJ)Z>$$. But this is nothing else 
that (2.1)!
\begin{pro}
Let $(M^{2n},g,J)$ be a complete, simply connected, NK-manifold. Then $M$ is a riemannian 
product $M_1 \times M_2$ where $M_1$ is a K\"ahler manifold and $M_2$ 
a strict NK-manifold.
\end{pro}
{\bf{Proof}} : \\
Set $E_1=Ker(r)$ and let $E_2$ be the orthogonal complement of $E_1$ in $TM$. By (2.1) 
both $E_1$ and $E_2$ are $\overline{\nabla}$-parallel. Since $\nabla_XJ$ vanishes 
whenever $X$ is in $E_1$ the distribution $E_1$ is in fact $\nabla$-parallel. Now, if 
$X$ is in $TM$ and $Y$ in $E_2$ we have $(\nabla_XJ)Y \in E_1^{\perp}=E_2$, hence 
$E_2$ is $\nabla$-parallel. It is now easy to conclude by a theorem of de Rham
$\blacksquare$ 
\begin{rema}
Proposition 2.1 was already proven in \cite{Gray1} under the assumption that the 
tensor $r$ is $\nabla$-parallel.
\end{rema}
\par
Therefore, we can restrict our attention to the class of strict NK-manifolds. 
\begin{pro}
Let $(M^{2n},g,J)$ a strict NK-manifold.  \\
(i) Suppose that $r$ has more than one eigenvalue. Then the canonical Hermitian connection has reduced holonomy. \\
(ii) If the tensor $r$ has exactly one eigenvalue then $M$ is a positive 
Einstein manifold. Furthermore, the first Chern class of $(M,J)$ vanishes.
\end{pro}
{\bf{Proof}} : \\
(i) If $\lambda_i>0, i=\overline{1,p}$ are the eigenvalues of $r$ we have a 
$\overline{\nabla}$-parallel decomposition 
\begin{nr} \hfill 
$ TM= \bigoplus \limits_{i=1}^{p} E_i \hfill $
\end{nr}
where $E_i$ is the eigenbundle corresponding to the eigenvalue $\lambda_i$. Hence, each factor is 
preserved by the holonomy group, which is thus reducible. \\
(ii) The proof can be found in \cite{Gray1}, page 242. Let us give it for the sake of completeness.
We recall the following formula : 
\begin{nr} \hfill 
$ \sum \limits_{i,j=1}^{2n}<re_i,e_j>(R(X,e_i,Y,e_j)-5R(X,e_i,JX,Je_j))=0\hfill $
\end{nr}
(see \cite{Gray1}) where $\{e_i\}_{i=\overline{1,2n}}$ is a local orthonormal 
frame field and  $X,Y$ are in $TM$. If $r=\lambda 1_{TM}, \lambda>0$ this formula becomes 
$Ric-5Ric^{\star}=0$ hence, $Ric=\frac{5\lambda}{4}$ as $Ric-Ric^{\star}=r$. 
The second assertion follows by the description of the first Chern class of $(M,J)$ 
given in \cite{Gray1}
$\blacksquare$ \\ \par
The first part of the theorem 1.1 follows now from the previous proposition. 
We will now compute the Ricci tensor of a NK-manifold and 
show that it is completely determined by the spectral decomposition of the tensor $r$. This 
computation will be equally used in section 5.
\begin{lema}We have, by respect to the decomposition (2.2) : \\
(i) $Ric(X,Y)=0$ if $X$ and $Y$ are vector fields belonging to $E_i$ and $E_j$ respectively, 
and $i \neq j$. \\
(ii) If $X,Y$ are vector fields in $E_i$ : 
$$ Ric(X,Y)=\frac{\lambda_i}{4}<X,Y>+\frac{1}{\lambda_i} \sum \limits_{s=1}^{p} 
\lambda_s <r^s(X),Y>$$
where the tensors $r^s : TM \to TM, 1 \le s \le p$ are defined by 
$<r^s(X),Y>=-Tr_{E_s}(\nabla_XJ)(\nabla_YJ)$ whenever $X,Y$ are in $TM$.
\end{lema}
{\bf{Proof}} : \\
(i) Let us denote by $\overline{R}$ the curvature tensor of the connexion $\overline{\nabla}$. 
We have (see \cite{Gray1}, page 237) :
\begin{nr} \hfill 
$ \begin{array}{lr}
\overline{R}(X,Y,Z,T)=R(X,Y,Z,T)-\frac{1}{2}<(\nabla_XJ)Y, 
(\nabla_ZJ)T>+ \\
+\frac{1}{4} \biggl [ <(\nabla_XJ)Z, (\nabla_YJ)T>-<(\nabla_XJ)T, (\nabla_YJ)Z> \biggr ].
\end{array} \hfill $
\end{nr}
Let $\{ e_k \}_{k=\overline{1,2n}}$ on orthonormal base of $TM$ which gives 
orthonormal bases in $E_s$ for $1 \le s \le p$. We get : 
$$Ric(X,Y)=\sum \limits_{s=1}^{p} \sum \limits_{e_k\in E_s}^{} R(X,e_k,Y,e_k).$$
If $s \neq j$ we have $\overline{R}(X,e_k,Y,e_k)=0$ hence $R(X,e_k,Y,e_k)=
\frac{1}{4}<((\nabla_{e_k}J))X,(\nabla_{e_k}J)Y>$ by (2.4). If $s=j$ then $s \neq i$ and 
as before we get $$R(X,e_k,Y,e_k)=R(Y,e_k,X,e_k)=
\frac{1}{4}<((\nabla_{e_k}J))X,(\nabla_{e_k}J)Y>.$$ It follows that 
$Ric(X,Y)=\frac{1}{4}<rX,Y>=0$.\\
(ii) Using (2.3) we obtain : 
$$ \sum \limits_{s=1}^{p} \lambda_s \biggl ( \sum 
\limits_{e_k \in E_s}^{} R(X,e_k,Y,e_k)-5R(X,e_k,JY,Je_k)\biggr )=0.$$
Reasoning as in the proof of (i), we get for $s \neq i$ that 
$$R(X,e_k,JY,Je_k)=-3R(X,e_k,Y,e_k)=-\frac{3}{4}<(\nabla_{e_k}J)X, (\nabla_{e_k}J)Y>.$$
It follows that 
$$ 4\sum \limits_{\stackrel{s=1}{s \neq i}}\lambda_s<r^sX,Y>+\lambda_i 
\biggl ( \sum \limits_{e_k \in E_s}^{} R(X,e_k,Y,e_k)-5R(X,e_k,JY,Je_k)\biggr )=0$$
and further 
$4\sum \limits_{\stackrel{s=1}{s \neq i}}(\lambda_s-\lambda_i)<r^sX,Y>+\lambda_i 
<(Ric-5Ric^{\star})X,Y>=0$. We conclude by using that $Ric-Ric^{\star}=r$ and 
$\sum \limits_{s=1}^{p}r^s=r$
$\blacksquare$ \\ \par
Note that by  definition the tensors $r^s, 1 \le s \le p$ are positive. Setting $\lambda=
min \{\lambda_i : 1 \le i \le p \}$ the proposition 2.1 obviously implies 
that $Ric \ge \lambda g$. This, together with Myer's theorem proves the second assertion 
of theorem 1.1. \par
Another result we will use in the next section is : 
\begin{lema}
The tensors $r^s, 1 \le s \le p$ are $\overline{\nabla}$-parallel.
\end{lema}
The proof is analogous to that of the $\overline{\nabla}$-parallelism of $r$
so it will be left to the reader. Thus, using the lemma 2.1 we obtain that :  
\begin{coro}
The Ricci tensor and the Ricci $\star$ tensor of a compact NK-manifold are 
$\overline{\nabla}$-parallel. 
\end{coro}
It follows that the scalar curvature and more, the $\star$-scalar curvature of $(M,g,J)$,
are strictly positive constants. The proof of the theorem 1.1 is now finished.
\section{Examples of NK manifolds}
Let us consider a 
Riemannian submersion with totally geodesic fibers 
$$F \hookrightarrow (M,g) \to N $$ and let $TM={\cal{V}} \oplus H$ be the corresponding splitting 
of $TM$. We will suppose that $M$ admits a complex structure $J$ compatible with 
$g$ and preserving ${\cal{V}}$ and $H$ such that $(M,g,J)$ is a K\"ahler manifold. 
Consider now the Riemannian metric on $M$ defined by 
$$ \hat{g}(X,Y)=\frac{1}{2}g(X,Y) \ \mbox{if} \ X,Y \in {\cal{V}},  
\hat{g}(X,Y)=g(X,Y) \ \mbox{for} \ X,Y \mbox{in} \ H.$$
The metric $\hat{g}$ admits a compatible almost complex structure 
$\hat{J}$ given by $\hat{J}_{\vert {\cal{V}}}=-J$ and $\hat{J}_{\vert H}=J$. This almost complex structure 
was introduced in \cite{Eells} for the case of twistor spaces over $4$-manifolds.
\begin{pro}
The manifold $(M,\hat{g},\hat{J})$ is nearly K\"ahler. The distributions 
${\cal{V}}$ and $H$ are parallel with respect to the canonical Hermitian 
connection of $(M,\hat{g},\hat{J})$ which thus has reduced holonomy.
\end{pro}
{\bf{Proof}} : 
Let $A : TM \times TM \to TM$ be the O'Neill tensor of the Riemannian submersion $(M,g)$. As 
$g$ is K\"ahler we must have $A_XJ=JA_X$ for all $X$ in $TM$. Using the relations 
between the Levi-Civita connections of $\hat{g}$ and $g$ given in \cite{Besse1} we obtain 
after a standard computation : 
$$ \begin{array}{lr}
(\hat{\nabla}_X \hat{J})V=-(\hat{\nabla}_V \hat{J})X=-A_X(JV) \\
(\hat{\nabla}_V \hat{J})W=0, \ (\hat{\nabla}_X \hat{J})Y=2A_X(JY)
\end{array} $$
for every $X,Y$ in ${\cal{V}}$ and $V,W$ in $H$. It is now straightforward to conclude
$\blacksquare$ 
\begin{coro}
The twistor space of a positive quaternionic-K\"ahler manifold of 
dimension $4k$ admits a canonical  NK structure with reducible holonomy, 
contained in $U(1) \times U(2k)$.
\end{coro}
{\bf{Proof}} : \\
We have only to recall \cite{Sal3} that such a twistor space is the total space of a Riemannian submersion with 
totally geodesic fibers of dimension $2$ and that it admits a compatible K\"ahler structure
$\blacksquare$

\section{Reducible NK manifolds}
In this section we consider strict NK-manifolds $(M^{2n},g,J)$ such that the holonomy of the 
canonical Hermitian connection is contained in $U(1) \times U(n-1)$. This leads to a
$\overline{\nabla}$-parallel decomposition 
of $TM$, orthogonal with respect to $g$ and stable by $J$ 
$$ TM=L \oplus E$$ 
with $L$ of rank two. Note that the torsion of $\overline{\nabla}$ vanishes on $L$ and 
$T(L,E) \subseteq E$. 
\begin{lema} We have : \\
(i) $\overline{R}(X,Y,V,JV)=
-2<(\nabla_VJ)^2X,JY>$ for every vector fields $X,Y$ on $E$ and $V$ on $L$. \\
(ii) $\overline{R}(X,V,V,JV)=0$ if $X$ belongs to $E$ and $V$ to $L$.
\end{lema}
{\bf{Proof}} : \\ 
(i) Using (2.4) we get 
$$ \overline{R}(X,Y,V,JV)=R(X,Y,V,JV)-\frac{1}{2}<(\nabla_VJ)^2X,JY>.$$
Now the first Bianchi identity gives $R(X,Y,V,JV)=-R(Y,V,X,JV)+R(X,V,Y,JV)$. As 
$E$ is $\overline{\nabla}$-parallel we must have $\overline{R}(Y,V,X,JV)=0$ so we find  
by (2.4) that $R(Y,V,X,JV)=\frac{3}{4}<(\nabla_VJ)^2X,JY>$. In the same way we have 
$R(X,V,Y,JV)=-\frac{3}{4}<(\nabla_VJ)^2X,JY>$ and the result follows easily. \\
(ii) Using  (2.4) twice we get  
$$\overline{R}(X,V,V,JV)=R(X,V,V,JV)=R(V,JV,X,V)=\overline{R}(V,JV,X,V)$$ and we conclude 
by the fact that $E$ is $\overline{\nabla}$-parallel $\blacksquare$ \\
\par
Let us denote by $\Omega$ the curvature form of the line bundle $L$. Then we have 
$$ \overline{R}(X,Y)V=\Omega(X,Y)JV$$
for $X,Y$ in $TM$ and $V$ in $L$. We denote by $\omega^L$ the restriction 
of the K\"ahler form $\omega$ to $L$. Let $F$ be the endomorphism of $TM$ defined by 
$<FX,Y>=-\frac{1}{2}Tr_L(\nabla_XJ) (\nabla_YJ)$ whenever $X,Y$ are in $TM$.
\begin{rema}
If $V$ is a local vector field on $L$ of norm $1$ we have $F=-(\nabla_VJ)^2$. Hence 
$F$ is symmetric and positive, with $[F,J]=0$. By lemma 2.2 $F$ is $\overline{\nabla}$-
parallel and it follows easily that $\nabla_VF=0$ for every vector field $V$ in $L$.
\end{rema}
If $q^E$ is the $2$-form on $E$ defined by $q^E(X,Y)=
<FX,JY>$ for $X,Y$ in $E$ we obtain by lemma 4.1 that : 
$$ \Omega=f \omega^L+2q^E $$
where $f$ is a smooth function on $M$.
\begin{lema} We have : \\
(i) $d \omega^L(X,V,JV)=dq^E(X,V,JV)=0$ if $V$ is in $L$ and $X$ in $E$.\\
(ii) $$ \begin{array}{lr} 
d \omega^L(V,X,Y)=-<(\nabla_VJ)X,Y> \\
d q^E(V,X,Y)=-2<F(\nabla_VJ)X,Y>
\end{array} $$ where $V,X,Y$ are vector fields belonging to $L$ resp. $E$.
\end{lema}
{\bf{Proof}} : \\
The proof of (i) is straightforward. We leave it to the reader and concentrate on (ii). We have 
$$d \omega^L(V,X,Y)= \nabla_V\omega^L(X,Y)-\nabla_X\omega^L(V,Y)+
\nabla_Y\omega^L(V,X).$$
The fact that $\omega^L$ vanishes as soon as we take a direction in $E$ gives us 
that \\
$\nabla_V\omega^L(X,Y)=0, \nabla_X\omega^L(V,Y)=-\omega^L(V,\nabla_XY)$ and 
$\nabla_Y\omega^L(V,X)=-\omega^L(V, \nabla_YX).$ The claimed formula for 
$d \omega^L(V,X,Y)$ follows using the fact that $\overline{\nabla}_XY$ and 
$\overline{\nabla}_YX$ belong to $E$. Next, we have 
$$ dq^E(V,X,Y)=(\nabla_Vq^E)(X,Y)-(\nabla_Xq^E)(V,Y)+(\nabla_Yq^E)(V,X).$$
The vanishing of $q^E$ on $L \times E$ implies that 
$$(\nabla_Vq^E)(X,Y)=
<(\nabla_VF)X,JY>+<FX,(\nabla_VJ)Y>=<FX,(\nabla_VJ)Y>$$ 
(see the remark 4.1) and 
$$(\nabla_Xq^E)(V,Y)=\frac{1}{2}<F(\nabla_VJ)X,Y>, \ 
(\nabla_Y q^E)(V,X)=\frac{1}{2}<F(\nabla_VJ)Y,X>.$$ We conclude by using the fact that 
$F$ commutes with $\nabla_VJ$ 
$\blacksquare$ \\
\par
Let $\omega^E$ be the restriction of the form $\omega$ to $E$. We can now have a 
complete description of the curvature form of our line bundle $L$ as follows. 
\begin{pro}(i) There exists a constant $k>0$ such that $F_{\vert E}=
\frac{k}{4}1_E$. Moreover, the 
curvature form of the line bundle $L$ is 
$$ \frac{k}{2}(-2\omega^L+\omega^E).$$
(ii) We have that $(\nabla_XJ)Y$ belongs to $L$ whenever $X,Y$ are in $E$.
\end{pro}
{\bf{Proof}} : \\
(i) As $\Omega$ is closed 
we get $fd\omega^L+df \wedge \omega^L=-2dq^E$. If $X$ resp. $V$ are vector fields in $E$ 
resp. $L$ it follows by lemma 4.2, (i) that $X.f=0$, hence 
$df_{\vert E}=0$. This implies that $[X,Y].f=0$ whenever 
$X,Y$ are vector fields in $E$ and further that $(\nabla_XJ)Y.f=0$ (here we used that 
$E$ is $\overline{\nabla}$-parallel and $\overline{\nabla}_XY-\overline{\nabla}_YX=
[X,Y]+(\nabla_XJ)JY$). But the map $u : E \times E \to L$ 
defined at $(v,w) \in E \times E$ as the orthogonal projection of 
$(\nabla_vJ)w$ on $L$  is surjective 
by the injectivity of $F_{\vert E_m}$. Hence $df$ vanishes on $L$ and thus $df=0$, that 
is $f$ is constant, equal to $c$.\par
Let now $X,Y$ resp. $V$ be vector fields in $E$ resp. $L$. As $d\Omega(V,X,Y)=0$ we get 
by lemma 4.2, (ii) 
$$ -c<(\nabla_VJ)X,Y>-4<F(\nabla_VJ)X,Y>=0.$$
We deduce that $(\nabla_VJ)(4F+c)=0$ and further $F(4F+c)=0$ on $E$. As the restriction of 
$F$ to $E$ is injective it follows that $F=\frac{-c}{4}id$ on $E$. We set $k=-c$. \\
(ii) Let $X,Y,Z$ be vector fields on $E$. As we obviously have $d\omega^L(X,Y,Z)=0$ it follows 
by (i) that $d \omega^E(X,Y,Z)=0$. A straightforward computation 
gives $(\nabla_X \omega^E)(Y,Z)=-<(\nabla_XJ)Y,Z>$ from which we deduce that 
$d \omega^E(X,Y,Z)=-<(\nabla_XJ)Y,Z>$ $\blacksquare$ \\ \par
\begin{coro}
(i) The tensor $r$ has exactly two eigenvalues : 
$\frac{k(n-1)}{2}$ resp. $k$ with eigenbundles $L$ 
resp. $E$. \\
(ii) The Ricci tensor of $(M,g)$ has exactly two eigenvalues : $\frac{k(n+7)}{8}$ and 
$\frac{k(n+2)}{4}$ with eigenbundles  $L$ resp. $E$.
\end{coro}
{\bf{Proof}} : \\
(i) The fact that $r_{\vert L}=\frac{k(n-1)}{2}$ follows easily by the fact that 
$F$ is constant on $E$. If 
$x$ is in $E$ let $v$ in $L$ be unitary, and $\{e_i \}_{1 \le i \le 2(n-1)}$ an orthogonal basis 
of $E$. Then we have  $<rx,x>=2 \Vert (\nabla_vJ)x\Vert^2+
\sum \limits_{i=1}^{2(n-1)} \Vert (\nabla_{e_i}J)x\Vert^2 $. As $(\nabla_{e_i}J)x$ belongs 
to $L$, the last sum equals $2 \Vert (\nabla_vJ)x\Vert^2$ and we use $F_{\vert E}=
\frac{k}{4}$. \\
(ii) Follows from lemma 2.1 and (i)
$\blacksquare$
\section{The twistor structure}
Let us define a new Riemannian metric on $M$, called $\overline{g}$, as follows : 
$$ \overline{g}(X,Y)=g(X,Y) \ \mbox{if} \ X,Y \in E, \ \overline{g}(X,Y)=2g(X,Y) \ \mbox{for} \ 
X,Y \ \mbox{in} \ L.$$
The {\it{reversing}} almost complex structure defined by $\overline{J}_{\vert L}=-J$ 
and $\overline{J}_{\vert E}=J$ is in fact integrable, the proof being identical to 
that given in six dimensions in \cite{Moro1}. The K\"ahler form of $(M, \overline{g}, \overline{J})$ is exactly 
$-2\omega^L+\omega^E$ and hence it is closed by proposition 4.1, (i). Thus, 
$(M, \overline{g}, \overline{J})$ is a K\"ahler manifold.
\begin{lema} $(M, \overline{g})$ is an Einstein manifold, with Einstein constant 
$\frac{n+1}{4}k$.
\end{lema}
{\bf{Proof}} : \\
This is a computation very similar to that of \cite{Besse1}, page 232, where the Ricci 
tensor of the canonical variation of a Riemannian submersion is computed. Let 
$\widetilde{\nabla}$ be the Levi-Civita connection of the metric $\overline{g}$. If 
$V$ resp. $X,Y$ are vector fields in $L$ resp. $E$ we have : 
$\widetilde{\nabla}_VX=\overline{\nabla}_VX, \widetilde{\nabla}_XV=\overline{\nabla}_XV-
(\nabla_XJ)JV$ and $\widetilde{\nabla}_XY=\nabla_XY$. Moreover, $\widetilde{\nabla}_VW=
\nabla_VW$  whenever $V,W$ are in $L$.
This follows from the definition
of the Levi-Civita connection and by the fact that the $\overline{\nabla}$-parallelism 
of $L$ and $E$ allows us to identify the projections on $L$ resp. $E$ of 
brackets of the type $[V,X]$ and $[X,Y]$.\par
Let $\widetilde{R}$ be the curvature tensor of $\widetilde{\nabla}$. Using the above formulas 
we get, after a standard computation : 
$$\begin{array}{lr}
\widetilde{R}(V,X,V,X)=<FX,X>\Vert V \Vert^2=\frac{k}{4}\Vert V \Vert^2 \Vert X \Vert^2 \\
\widetilde{R}(X,Y,X,Y)=\overline{R}(X,Y,X,Y)-\frac{1}{2}\Vert (\nabla_XJ)Y\Vert^2=
R(X,Y,X,Y)-\frac{3}{4}\Vert (\nabla_XJ)Y\Vert^2
\end{array} $$
by (2.4). The result follows now by corollary 4.1
$\blacksquare$ \\ \par
Thus, 
$(M, \overline{g}, \overline{J})$ is a K\"ahler-Einstein manifold, which is also Fano. Moreover, 
the distribution $E$ defines a complex contact structure on the 
complex manifold $(M,\overline{J})$ as it is 
$\overline{J}$-holomorphic and the map $(X,Y) \in E \times E \to 
(\nabla_XJ)Y$ which gives the Frobenius obstruction 
is everywhere non-degenerate. By a result of 
LeBrun (see \cite{Lebrun1}) $(M, 
\overline{g})$ is the twistor space of a positive 
quaternionic-K\"ahler manifold. Moreover, from the construction of the metric $\overline{g}$ we deduce 
that $(M,g)$ is is the twistor space of a positive 
quaternionic-K\"ahler manifold endowed with its canonical NK structure. This proves theorem 1.2.
\begin{rema}
If $M$ is of dimension $6$, it has constant type and proposition 4.1 is automatically 
satisfied. Corollary 4.1 follows by the fact that every $6$-dimensional NK 
manifold is Einstein \cite{Gray1}. Thus all we need to prove the theorem 1.1 in 
this case is lemma 5.1.
\end{rema}\par
Let us prove now the corollary 1.1. It is well known (see \cite{Gray1}) that in $10$-dimensions 
the eigenvalues of $r$ are $4(\alpha^2+\beta^2)$ with multiplicity $2$, $4\alpha^2$ and 
$4\beta^2$ each of multiplicity $4$, where $\alpha \ge \beta \ge 0$. If $\beta=0$ then it follows 
by \cite{Gray1} that the universal cover of $M$ is a Riemannian product as stated. If $\beta>0$ then $M$ is 
strict and we apply theorem 1.1.

$\\$
\begin{flushright}
Paul-Andi Nagy \\
Institut de Math\'ematiques, Universit\'e de Neuch\^atel, \\
rue E. Argand 11, 2007 Neuch\^atel, Switzerland 
\\ e-mail : Paul.Nagy@unine.ch
\end{flushright}

\normalsize
\begin{thebibliography}{99}
\bibitem{Al1}
B.ALEXANDROV, G.GRANTCHAROV, S.IVANOV, \textit{Curvature properties of twistor 
spaces of quaternionic K\"ahler manifolds}, J.Geom. {\bf{62}} (1998), 1-12.
\bibitem{Moro1}
F. BELGUN, A. MOROIANU, \textit{Nearly K\"ahler $6$-manifolds with reduced holonomy}, 
Ann. Glob.  An. Geom. {\bf{19}} (2001), 307-319.
\bibitem{Besse1}
A. BESSE, \textit{Einstein manifolds}, Springer-Verlag, New York, 1987. 
\bibitem{Eells}
J. EELS, S. SALAMON, \textit{Constructions twistorielles des applications harmoniques}, C. R. Acad. Sc. Paris, 
{\bf{296}} (1983), 685-687.
\bibitem{Friedrich2} 
Th. FRIEDRICH, S. IVANOV, \textit{Parallel spinors and connections with skew-symmetric 
torsion in string theory}, math. DG/0102142 (2001).
\bibitem{Gray5}
A. GRAY, \textit{Nearly K\"ahler manifolds}, J. Diff. Geom. {\bf{4}} (1970), 283-309
\bibitem{Gray2}
A. GRAY, \textit{Weak holonomy groups}, Math. Z, {\bf{125}} (1971), 290-300.
\bibitem{Gray4}
A. GRAY, \textit{Riemannian manifolds with geodesic symmetries of order 3}, J. Diff. 
Geometry {\bf{7}} (1972), 343-369.
\bibitem{Gray1} 
A. GRAY, \textit{The structure of nearly K\"ahler manifolds}, Math. Ann. 
{\bf{223}} (1976), 233-248.
\bibitem{Gray3}
A. GRAY, L. M. HERVELLA, \textit{The sixteen classes of almost Hermitian manifolds and their 
linear invariants}, Ann. Mat. Pura Appl. {\bf{123}} (1980), 35-58.
\bibitem{Grun1}
R. GRUNEWALD, \textit{Six-Dimensional Riemannian manifolds with real Killing 
spinors}, Ann. Global Anal. Geom. {\bf{8}} (1990), 43-59.
\bibitem{Sal2}
C. LEBRUN, S. SALAMON, \textit{Strong rigidity of positive quaternion-K\"ahler 
manifolds}, Invent. Math. {\bf{118}} (1994), 109-132. 
\bibitem{Lebrun1}
C. LEBRUN, \textit{Fano manifolds, contact structures and quaternionic geometry}, 
International J. Math. {\bf{6}} (1995), 419-437.
\bibitem{Mu}
O. MUSKAROV, \textit{Structures presque Hermitiennes sur les espaces twistoriels et leurs types}, 
C. R. Sc. Paris Ser. I Math. {\bf{305}} (1987), 365-398.
\bibitem{Sal1}
Y. S. POON, S.SALAMON \textit{Eight-dimensional quaternion K\"ahler manifolds 
with positive scalar curvature}, J.Differential. Geom. {\bf{33}} (1990), 363-378.
\bibitem{Sal3}
S. SALAMON, \textit{Quaternionic K\"ahler manifolds}, Invent. Math. {\bf{67}} (1982), 143-171.
\bibitem{Swann1}
A. SWANN, \textit{Weakening holonomy}, preprint ESI no. 816, (2000)
\end{thebibliography}
\end{document}